# WHAT IS THE PROBABILITY OF INTERSECTING THE SET OF BROWNIAN DOUBLE POINTS?

By Robin Pemantle[1] and Yuval Peres[2]

*University of Pennsylvania and University of California*

We give potential theoretic estimates for the probability that a set $A$ contains a double point of planar Brownian motion run for unit time. Unlike the probability for $A$ to intersect the range of a Markov process, this cannot be estimated by a capacity of the set $A$. Instead, we introduce the notion of a capacity with respect to two gauge functions simultaneously. We also give a polar decomposition of $A$ into a set that never intersects the set of Brownian double points and a set for which intersection with the set of Brownian double points is the same as intersection with the Brownian path.

**1. Introduction.** Let $A$ be a compact subset of the $\frac{1}{3}$-unit disk in the plane. For fifty years it has been known that $A$ intersects the path of a Brownian motion with positive probability if and only if $A$ has positive Newtonian capacity. In fact, the Newtonian (logarithmic) capacity gives an estimate, up to a constant factor, the probability that $A$ is hit by a Brownian motion started, say, from the point $(1,0)$ and run for a fixed time. The estimate is of course stronger than the dichotomous result, and moreover, it turns out to be important when examining properties of intersections with random sets; see, for example, the simple Cantor-type random fractal shown in Peres [5] to be "intersection-equivalent" to the Brownian motion; see also the remark after Theorem 2.3.

Similar results are known for much more general Markov processes. Let $G(x,y)$ denote the Green function for a transient Markov process. The capacity, $\text{Cap}_K(A)$ of a set $A$, with respect to a kernel $K$ is defined to be the

Received November 2003; revised May 2006.

[1]Supported in part by NSF Grant DMS-01-03635.

[2]Supported in part by NSF Grant DMS-98-03597.

reciprocal of the infimum of energies

$$\mathcal{E}_K(\mu) := \int\int K(x,y)\,d\mu(x)\,d\mu(y)$$

as $\mu$ ranges over probability measures supported on $A$. In a wide variety of cases it is known that the range of the process intersects $A$ with positive probability if and only if $A$ has positive capacity with respect to the Green kernel. The same is true of any of a number of related kernels, and choosing the Martin kernel $M(x,y) = G(x,y)/G(\rho,y)$ with respect to any starting point $\rho$ (see, e.g., Benjamini, Pemantle and Peres [1]) leads to the estimate

$$\tfrac{1}{2}\operatorname{Cap}_M(A) \leq \mathbf{P}_\rho(\text{the process intersects } A) \leq \operatorname{Cap}_M(A).$$

We are chiefly interested in the set $\mathcal{D}$ of double points of a planar Brownian motion. We work on a probability space $(\Omega, \{\mathcal{F}_t\}, \mathbf{P})$ on which are defined two independent Brownian motions, $B_t$ and $\tilde{B}_t$, both started from the point $\rho := (1,0)$. The notation $\mathbf{P}_x$ (or $\mathbf{P}_{x,y}$) will be used when a different starting point (or points) is required. Let $\tau_* = \inf\{t : |B_t| = 3\}$ be the exit time of $B_t$ from the disk $\{|x| \leq 3\}$. Formally, then,

$$\mathcal{D} := \{x : B_r = B_s = x \text{ for some } 0 < r < s < \tau_*\}.$$

The choice to start at $\rho$, stop at $\tau_*$, and choose sets inside the $\tfrac{1}{3}$-unit disk are conveniences that make the Martin and Green kernel both comparable to $|\log|x-y||$.

The random set $\mathcal{D}$ is not the range of any Markov process, but we may still ask about the probability for the random set $\mathcal{D}$ to intersect a fixed set $A$. A closely related random set to $\mathcal{D}$ is the intersection of two independent Brownian motions, denoted here by

$$\mathcal{I} := \{x : B_r = \tilde{B}_s = x \text{ for some } 0 < r < \tau_*, 0 < s < \tilde{\tau}_*\},$$

where $\tilde{\tau}_* = \inf\{t : |\tilde{B}_t| \geq 3\}$. Fitzsimmons and Salisbury [3] showed, for a subset $A$ of the $\tfrac{1}{3}$-unit disk, that $\mathbf{P}(\mathcal{I} \cap A \neq \varnothing)$ may be estimated up to a constant factor by $\operatorname{Cap}_L(A)$ where $L(x,y) = (\log|x-y|)^2$. In general, they show that taking intersections of random sets multiplies the kernels in the capacity tests; see also Salisbury [6] and Peres [5]. The set $\mathcal{D}$ may be written as a countable union of the sets of $\varepsilon$-separated double points (we use a time separation of $\varepsilon^2$ so that $\varepsilon$ may be thought of as a small spatial unit):

$$\mathcal{D}_\varepsilon := \{x : B_r = B_s = x \text{ for some } 0 < r < r + \varepsilon^2 \leq s < \tau_*\}.$$

It is not hard to see that each random set $\mathcal{D}_\varepsilon$ behaves similarly to the set $\mathcal{I}$, but with an increasingly poor constant. In other words,

$$c_\varepsilon \operatorname{Cap}_L(A) \leq \mathbf{P}(\mathcal{D}_\varepsilon \cap A \neq \varnothing) \leq C_\varepsilon \operatorname{Cap}_L(A),$$

but the constant $C_\varepsilon$ goes to infinity as $\varepsilon$ goes to zero. Since the property of having zero capacity is closed under countable unions, we again have the dichotomous criterion

$$(1.1) \qquad \mathbf{P}(\mathcal{D} \cap A \neq \varnothing) = 0 \Leftrightarrow \mathrm{Cap}_L(A) = 0$$

for $L(x,y) = (\log|x-y|)^2$. No estimate follows, however.

An example helps to explain this shortcoming. Fix an $\alpha \in (1/2, 1)$ and let $A_n$ be nested subsets of the line segment $A_0 := [-1/2, 1/2] \times \{0\}$ such that $A_n$ is made of $2^n$ intervals of length $2^{-2^{\alpha n}}$, with each of the $2^n$ intervals of $A_n$ containing exactly two intervals of $A_{n+1}$ situated at the opposite ends of the interval of $A_n$. The intersection, denoted $A$, is a Cantor set for which, if $K(x,y) = |\log|x-y||$ and $L(x,y) = \log^2|x-y|$, then

$$\mathrm{Cap}_K(A) > 0 = \mathrm{Cap}_L(A).$$

For each set $A_n$, a Brownian motion that hits the set will immediately after have a double point in the set. Thus,

$$\mathbf{P}(\mathcal{D} \cap A_n \neq \varnothing) := p_n,$$

where $p_n$ decreases as $n \to \infty$ to a positive number, estimated by $\mathrm{Cap}_K(A)$. On the other hand, since $\mathrm{Cap}_L(A) = 0$, we know that $\mathcal{D}$ is almost surely disjoint from $A$.

From this we see that the probability of $A$ intersecting $\mathcal{D}$ is not continuous as $A$ decreases to a given compact set, and therefore, that this probability cannot be uniformly estimated by $\mathrm{Cap}_K$ for any $K$, since $\mathrm{Cap}_K$ is a Choquet capacity, and must be continuous with respect to this kind of limit. On the other hand, since the probability that $A_n$ intersects $\mathcal{D}_\varepsilon$ is estimated by the Choquet capacity $\mathrm{Cap}_L(A_n)$ which goes to zero as $n \to \infty$, we see that these estimates are indeed getting worse and worse as $n \to \infty$ for fixed $\varepsilon$, and are only good when $\varepsilon \to 0$ as some function of $n$.

We remark that such behavior is possible only because $\mathcal{D}$ is not a closed set. Indeed, if $X$ is a random closed set and $\{Y_n\}$ are closed sets decreasing to $Y$, then the events $\{X \cap Y_n \neq \varnothing\}$ decrease to the event $\{X \cap Y \neq \varnothing\}$, whence

$$(1.2) \qquad \mathbf{P}(X \cap Y_n \neq \varnothing) \downarrow \mathbf{P}(X \cap Y \neq \varnothing).$$

The goal of this note is to provide a useful estimate for $\mathbf{P}(\mathcal{D} \cap A \neq \varnothing)$. We have just seen that it cannot be of the form $\mathrm{Cap}_K$ for some kernel, $K$. Instead, we must introduce the notion of a capacity with respect to two different kernels, which we denote $\mathrm{Cap}_{f \to g}$. We go about this two different ways. The first approach is to show that $\mathrm{Cap}_{f \to g}$ gives estimates on probabilities of intersection with $\mathcal{D}_\varepsilon$ which are uniform in $\varepsilon$ and thus allow passage to the limit. This relies on the result of Fitzsimmons and Salisbury (or Peres), so is less self-contained, but yields as a by-product the estimates for $\varepsilon > 0$ which may be considered interesting in themselves. The second is a softer and more elementary argument, which produces a sort of polar decomposition of the set $A$ but is less useful for computing. Section 2 states our results, Section 3 contains proofs of the estimates and Section 4 contains the proof of the decomposition result.

**2. Results.** Since Brownian motion is isotropic, we will restrict attention to kernels $K(x,y) = f(|x-y|)$ that depend only on $|x-y|$. When $K$ has this form, we write $\mathcal{E}_f$ and $\text{Cap}_f$ instead of $\mathcal{E}_K$ and $\text{Cap}_K$. Let $f$ and $g$ be functions from $\mathbb{R}^+$ to $\mathbb{R}^+$ going to infinity at zero, with $f \leq g$. Let $h_\varepsilon$ denote the function on $\mathbb{R}^+$ defined by

$$h_\varepsilon(x) = \begin{cases} f(x), & \text{if } x \geq \varepsilon, \\ g(x) \cdot \dfrac{f(\varepsilon)}{g(\varepsilon)}, & \text{if } x < \varepsilon. \end{cases}$$

Let $\text{Cap}_\varepsilon$ denote $\text{Cap}_{h_\varepsilon}$. The following result defines the hybrid capacity $\text{Cap}_{f \to g}$ as a limit and also characterizes it as "$\text{Cap}_f$ measured only at places where $\text{Cap}_g$ is positive."

PROPOSITION 2.1. *The limit $\lim_{\varepsilon \to 0} \text{Cap}_\varepsilon(A)$ exists. Denoting this limit by $\text{Cap}_{f \to g}(A)$, we have*

(2.1) $\quad \text{Cap}_{f \to g}(A) = [\inf\{\mathcal{E}_f(\mu) : \mathcal{E}_g(\mu) < \infty \text{ and } \mu(A) = 1\}]^{-1}.$

PROOF. If $\text{Cap}_g(A) = 0$, then both sides of (2.1) are clearly zero, so assume that $\text{Cap}_g(A) > 0$. For each $\varepsilon$, let $\mu_\varepsilon$ be a probability measure on $A$ that minimizes $\mathcal{E}_{h_\varepsilon}$, so that $\text{Cap}_{h_\varepsilon}(A) = \mathcal{E}_{h_\varepsilon}(\mu_\varepsilon)$. Since $f \leq h_\varepsilon$ for all $\varepsilon$, we have

$$\mathcal{E}_{h_\varepsilon}(\mu_\varepsilon) \geq \mathcal{E}_f(\mu_\varepsilon).$$

Observe that each $\mu_\varepsilon$ has finite $g$-energy and take the infimum on the left-hand side and the supremum on the right-hand side, then invert, to see that

$$\sup_\varepsilon \text{Cap}_\varepsilon(A) \leq [\inf\{\mathcal{E}_f(\mu) : \mathcal{E}_g(\mu) < \infty \text{ and } \mu(A) = 1\}]^{-1}.$$

On the other hand, if $\mu$ is any measure of finite $g$-energy, then by choice of $\mu_\varepsilon$, we know that

$$\mathcal{E}_{h_\varepsilon}(\mu_\varepsilon) \leq \mathcal{E}_{h_\varepsilon}(\mu).$$

As $\varepsilon \to 0$, dominated convergence shows that the right-hand side of this converges to $\mathcal{E}_f(\mu)$, and hence, that

$$\liminf_{\varepsilon \to 0} \text{Cap}_\varepsilon(A) \geq [\inf\{\mathcal{E}_f(\mu) : \mathcal{E}_g(\mu) < \infty \text{ and } \mu(A) = 1\}]^{-1},$$

which finishes the proof. $\square$

REMARK. The infimum in (2.1) need not be achieved. For example, if $A$ is a small disk, $f(x) = |\log x|$, and $g(x) = x^{-\alpha}$ for any $\alpha \in [1,2)$, then the infimum of logarithmic energies of probability measures on $A$ is equal to the log-energy of normalized one-dimensional Lebesgue measure on the boundary of the disk, and is strictly less than the logarithmic energy of any measure of finite $g$-energy.

This proposition is our only general result on hybrid capacities. For the remainder of the paper, $f$ will always be $|\log \varepsilon|$ and $g$ will always be $\log^2 \varepsilon$, so the notation $h_\varepsilon$ will be unambiguous. [We have also found the notation easier to read if we use $\log |x-y|/\log \varepsilon$ rather that $|\log|x-y||/|\log \varepsilon|$ or $\log(1/|x-y|)/\log(1/\varepsilon)$ whenever the signs cancel.] Our main interest in $\text{Cap}_\varepsilon$ is that it gives the estimate on the probability of an intersection with $\mathcal{D}_\varepsilon$.

THEOREM 2.2 (Estimates for intersecting $\mathcal{D}_\varepsilon$). *Let $f(x) = |\log x|$ and $g(x) = \log^2 x$. There are constants $c$ and $C$ such that, for any $\varepsilon > 0$ and any closed subset $A$ of disk $\{x: |x| \leq 1/3\}$,*

$$c\,\text{Cap}_\varepsilon(A) \leq \mathbf{P}(\mathcal{D}_\varepsilon \cap A \neq \varnothing) \leq C\,\text{Cap}_\varepsilon(A).$$

Since $\text{Cap}_\varepsilon \uparrow \text{Cap}_{f \to g}$ and $\mathcal{D} = \bigcup \mathcal{D}_\varepsilon$, our first main result follows as an immediate corollary.

THEOREM 2.3 (Two-gauge capacity estimate). *For the same constants $c$ and $C$, and the same $f$ and $g$,*

$$c\,\text{Cap}_{f \to g}(A) \leq \mathbf{P}(\mathcal{D} \cap A \neq \varnothing) \leq C\,\text{Cap}_{f \to g}(A).$$

REMARK. Suppose the set $A$ is a bi-Hölder image of some set $S$ for which the intersection probabilities with $\mathcal{D}$ are known. Since the logarithm of the distance between two points in a small disk changes by a bounded factor under such a map, the Newtonian and $\log^2$ capacities change only by a bounded factor, so the probability of $A$ intersecting $\mathcal{D}$ is estimated by the probability of $S$ intersecting $\mathcal{D}$. This is more than can be concluded from the dichotomy (1.1).

The characterization of $\text{Cap}_{f \to g}$ in Proposition 2.1 suggests an explanation for the two-gauge capacity result. The probability of intersection with $\mathcal{D}$ is estimated by $\text{Cap}_{\log}$ "at places of finite $\log^2$-energy," so perhaps the operative mechanism is that one must eliminate certain "thin" places that can never contain Brownian double points, leaving a "core set," such that if and when Brownian motion hits the core set, immediately there will be a Brownian double point in the core set. This turns out to be true.

THEOREM 2.4 (Polar decomposition). *Any compact subset $A$ of the plane not containing $(1,0)$ may be written as a union $A = A_1 \cup A_2$, such that (1) the set $A_1$ is almost surely disjoint from $\mathcal{D}$, and (2), on the event that the hitting time $\tau_2$ of $A_2$ is finite, then for any $\varepsilon > 0$, with probability 1, Brownian motion stopped at time $\tau_2 + \varepsilon$ has a double point in $A_2$.*

It follows from this that

$$\mathbf{P}(\mathcal{D} \cap A \neq \varnothing) = \mathbf{P}(\text{Brownian motion hits } A_2),$$

which is estimated up to a constant factor by $\text{Cap}_{\log}(A_2)$, and, in fact, is equal to the Martin capacity of $A_2$. Thus, this decomposition is in some ways stronger than Theorem 2.3; it is, in principle, less useful for computation because $A_2$ must first be computed, though, in practice, usually $A_2 = A$ or is empty. We remark that $\text{Cap}_{\log}(A_2)$ is a different estimate from $\text{Cap}_{\log \to \log^2}(A)$, if harmonic measure on $A_2$ has infinite $\log^2$-energy.

**3. Proof of estimates for intersecting $\mathcal{D}$.** Fix $\varepsilon \in (0, 1/3)$ and any $\delta < \varepsilon/2$. Let $x$ and $y$ be points in the quarter unit disk with $|x - y| > 3\delta$ and denote by $D_x$ and $D_y$ the balls of radius $\delta$ centered at $x$ and $y$, respectively. The key estimates for applying potential theoretic methods are the first and second moment estimates, as given in the following lemma. The notation $\asymp$ denotes equivalence up to a constant multiple.

LEMMA 3.1. *Let $H(A) = H(A, \varepsilon)$ denote the event $\{\mathcal{D}_\varepsilon \cap A \neq \varnothing\}$:*

$$(3.1) \qquad \mathbf{P}(H(D_x)) \asymp \frac{|\log \varepsilon|}{\log^2 \delta}.$$

*Letting $\mathbf{P}_\xi$ denote probabilities with respect to a Brownian motion started at the point $\xi \notin D_x$, we have, in general,*

$$(3.2) \qquad \mathbf{P}_\xi(H(D_x)) \asymp \frac{\log \varepsilon \log |\xi - x|}{\log^2 \delta}.$$

*The probabilities for double points simultaneously occurring in two balls are given as follows. When $|x - y| \geq \varepsilon$,*

$$(3.3) \qquad \mathbf{P}(H(D_x) \cap H(D_y)) \asymp \frac{|\log |x - y|| \cdot \log^2 \varepsilon}{\log^4 \delta}.$$

*When $|x - y| < \varepsilon$,*

$$(3.4) \qquad \mathbf{P}(H(D_x) \cap H(D_y)) \asymp \frac{|\log \varepsilon| \cdot \log^2 |x - y|}{\log^4 \delta}.$$

PROOF. Let $\tau$ be the hitting time on $D_x$. For $H(D_x)$ to occur, it is necessary that $\tau < \infty$ and that the Brownian motion hit $D_x$ after time $\tau + \varepsilon^2$. Denoting this event by $G$, use the Markov property at time $\tau$ and $\tau + \varepsilon^2$ and average over the position at time $\tau + \varepsilon^2$ to see that

$$\mathbf{P}(G) \asymp \frac{1}{|\log \delta|} \frac{\log \varepsilon}{\log \delta}.$$

On the other hand, conditioning on the position at time $\tau$ and at the return time to $D_x$, it is easy to bound $\mathbf{P}(H(D_x) \,|\, G)$ away from zero, since this is the probability that a Brownian path and a Brownian bridge, each started on the boundary of a ball of radius $\delta$ and run for time greater than $\delta^2$, intersect inside the ball. This establishes (3.1). When starting at a point $\xi$ near $x$ instead of at the point $(1,0)$, the probability of the event $\{\tau < \infty\}$ is $\log |\xi - x|/\log \delta$ rather than $1/|\log \delta|$, which gives the estimate in (3.2).

To establish the other two estimates, we consider possible sequences of visits, two to each ball, with the correct time separations. Let $H_1(x,y)$ denote the event that there exist times $0 < r < r + \varepsilon^2 \le s < t < t + \varepsilon^2 \le u < \tau_*$ such that $B_r \in D_x$, $B_s \in D_x$, $B_t \in D_y$ and $B_u \in D_y$. Let $H_2(x,y)$ denote the event that there exist times $0 < r < s < t < u < \tau_*$ such that $r + \varepsilon^2 \le t, s + \varepsilon^2 \le u$, $B_r \in D_x$, $B_s \in D_y$, $B_t \in D_x$ and $B_u \in D_y$. Let $H_3(x,y)$ denote the event that there exist times $0 < r < s < s + \varepsilon^2 \le t < u < \tau_*$ such that $B_r \in D_x$, $B_s \in D_y$, $B_t \in D_y$ and $B_u \in D_x$. The estimate

(3.5)    $\mathbf{P}(H(D_x) \cap H(D_y)) \asymp \mathbf{P}(H_1(x,y)) + \mathbf{P}(H_2(x,y)) + \mathbf{P}(H_3(x,y))$

follows from the same considerations: that for $j = 1, 2, 3$, $\mathbf{P}(H(D_x) \cap H(D_y) \,|\, H_j(x,y))$ is bounded away from zero; that the same holds when $x$ and $y$ are switched; that $\mathbf{P}(H_j(x,y)) \asymp \mathbf{P}(H_j(y,x))$; and that $H(D_x) \cap H(D_y)$ entails either $H_j(x,y)$ or $H_j(y,x)$ for some $j$. The estimates (3.3) and (3.4) will then follow from

(3.6) $$\mathbf{P}(H_1(x,y)) \asymp \frac{|\log |x-y|| \cdot \log^2 \varepsilon}{\log^4 \delta},$$

(3.7) $$\mathbf{P}(H_3(x,y)) \asymp \frac{\log^2 |x-y| \cdot |\log \varepsilon|}{\log^4 \delta},$$

(3.8) $$\mathbf{P}(H_2(x,y)) = O(\mathbf{P}(H_1(x,y)) + \mathbf{P}(H_3(x,y))).$$

The Markov property gives a direct estimate of $\mathbf{P}(H_1(x,y))$. In particular, we may take $r$ to be the hitting time of $D_x$, $s$ to be the next time after $r + \varepsilon^2$ that $D_x$ is hit, and so forth. The probability of hitting $D_x$ is $\asymp 1/|\log \delta|$. Given that $B_r \in D_x$, the probability that $B_s \in D_x$ for some $s \ge r + \varepsilon^2$ is $\asymp |\log \varepsilon|/|\log \delta|$. Given that, the probability of subsequently hitting $D_y$ is $\asymp |\log |x-y||/|\log \delta|$, and given such a hit at time $t$, the probability of $B_u \in D_y$ for some $u \ge t + \varepsilon^2$ is $\asymp |\log \varepsilon|/|\log \delta|$. Multiplying these together produces the estimate (3.6). Similarly, $\mathbf{P}(H_3(x,y))$ is the product of four factors, respectively comparable to $1/|\log \delta|$, $\log|x-y|/\log \delta$, $\log \varepsilon/\log \delta$ and $\log|x-y|/\log \delta$, proving (3.7).

In the case $|x-y| \geq \varepsilon$, the bound $\mathbf{P}(H_2(x,y)) = O(\frac{\log^3 |x-y|}{\log^4 \delta})$ is good enough to imply (3.8) and follows in the same manner from the Markov property at the hitting time of $D_x$, the next hit of $D_y$, the next hit of $D_x$ and the next hit on $D_y$. In the case $|x-y| \leq \varepsilon$, define an event $H_2' \subseteq H_2$ by additionally requiring $t \geq s + \varepsilon^2/2$. Let $H_2'' = H_2 \setminus H_2'$. The Markov property gives

$$(3.9) \qquad \mathbf{P}(H_2') = O\left(\frac{1}{|\log \delta|} \frac{\log|x-y|}{\log \delta} \frac{\log \varepsilon}{\log \delta} \frac{\log|x-y|}{\log \delta}\right) \asymp \mathbf{P}(H_3).$$

Finally, to estimate $\mathbf{P}(H_2'')$, observe that $H_2''$ entails both $s \geq r + \varepsilon^2/2$ and $u \geq t + \varepsilon^2/2$. The Markov property then gives

$$(3.10) \qquad \mathbf{P}(H_2'') = O\left(\frac{1}{|\log \delta|} \frac{\log \varepsilon}{\log \delta} \frac{\log|x-y|}{\log \delta} \frac{\log \varepsilon}{\log \delta}\right) \asymp \mathbf{P}(H_1)$$

and adding (3.9) to (3.10) establishes (3.8) and the lemma. $\square$

3.1. *Proof of the first inequality of Theorem* 2.2. The first inequality follows from Lemma 3.1 by standard methods. We give the details, since it is a little unusual to discretize space in only part of the argument (composing the set $A$ of lattice squares, but not discretizing the double point process itself). For the remainder of the argument, $\varepsilon$ and $A$ are fixed.

Let $\mu$ be any probability measure on $A$; we need to show that $\mathbf{P}(H(A)) \geq c\mathcal{E}_{h_\varepsilon}(\mu)^{-1}$. The closed set $A$ may be written as a decreasing intersection over finer and finer grids of finite unions of lattice squares. According to (1.2), we may therefore assume that $A$ is a finite union of lattice squares of width $\delta < \varepsilon$. Index the rows and columns of the grid, and let $\mathcal{B}$ denote the subcollection of squares where both coordinates are even. Let $\mathcal{B}'$ denote the collection of inscribed disks of $\mathcal{B}$. Then some translation $\mathcal{B}''$ of $\mathcal{B}'$ has $\mu$-measure at least $1/8$ (since space may be covered by 8 translates of the set of disks centered at points with both coordinates even). Define a random variable

$$X := \sum_{S \in \mathcal{B}''} \frac{\log^2 \delta}{|\log \varepsilon|} \mu(S) \mathbf{1}_{H(S)}.$$

By the first estimate in Lemma 3.1, the expectation of each $(\log^2 \delta/|\log \varepsilon|)\mathbf{1}_{H(S)}$ is bounded above and below by some constants $c_1$ and $c_2$. Thus, $c_1/8 \leq \mathbf{E}X \leq c_2$.

The second moment of $X$ is computed as

$$(3.11) \qquad \mathbf{E}X^2 = \frac{\log^4 \delta}{\log^2 \varepsilon} \sum_{S,T \in \mathcal{B}''} \mu(S)\mu(T) \mathbf{E}\mathbf{1}_{H(S) \cap H(T)}.$$

By estimates (3.3) and (3.4) of Lemma 3.1, when $S \neq T$,

$$\mathbf{E}\frac{\log^4 \delta}{\log^2 \varepsilon}\mathbf{1}_{H(S) \cap H(T)} \tag{3.12}$$

is bounded between constant multiples of $h_\varepsilon(|x - y|)$, where $x$ and $y$ are the centers of $S$ and $T$. Since $S$ and $T$ are separated by $\delta$, this is bounded between $c_3 h_\varepsilon(|x - y|)$ and $c_4 h_\varepsilon(|x - y|)$ for any $x \in S$ and $y \in T$. Thus, letting $U$ denote the union of $\mathcal{B}''$, the sum of the off-diagonal terms of (3.11) is estimated by

$$c_3 \int h_\varepsilon(x, y)\mathbf{1}_{|x-y|>\delta}\, d\mu(x)\, d\mu(y) \leq \frac{\log^4 \delta}{\log^2 \varepsilon} \sum_{S,T \in \mathcal{B}''} \mu(S)\mu(T)\mathbf{1}_{H(S) \cap H(T)}\mathbf{1}_{S \neq T}$$

$$\leq c_4 \int h_\varepsilon(x, y)\mathbf{1}_{|x-y|>\delta}\, d\mu(x)\, d\mu(y).$$

The diagonal terms sum to exactly $\mathbf{E}X$, so we see that

$$\mathbf{E}X^2 \leq \mathbf{E}X + c_4 \mathcal{E}_{h_\varepsilon}(\mu).$$

The second moment inequality $\mathbf{P}(X > 0) \geq (\mathbf{E}X)^2 / \mathbf{E}X^2$ now implies that

$$\mathbf{P}(X > 0) \geq \frac{c_1^2}{64(c_2 + c_4 \mathcal{E}_{h_\varepsilon}(\mu))}.$$

Since $X > 0$ implies the existence of an $\varepsilon$-separated double point in $A$, we have proved the first inequality with $c = c_1^3/(64(8c_2 + c_1 c_4))$.

3.2. *Proof of the second inequality of Theorem* 2.2. The following two propositions represent most of the work in finishing the proof of Theorem 2.2.

PROPOSITION 3.2. *If $A$ has diameter at most $\varepsilon$, then*

$$\mathbf{P}(\mathcal{D}_\varepsilon \cap A \neq \varnothing) \asymp |\log \varepsilon| \mathbf{P}(\mathcal{I} \cap A \neq \varnothing).$$

PROPOSITION 3.3 (Capacity criterion for $\mathcal{I}$). *For any $A$ in the $\frac{1}{3}$-unit disk,*

$$\mathbf{P}(\mathcal{I} \cap A \neq \varnothing) \asymp \mathrm{Cap}_{\log^2}(A).$$

The second of these two propositions is proved in Peres [5] but also follows from the methods of Fitzsimmons and Salisbury [3] if one upgrades to a quantitative estimate by observing that the Green kernel is comparable to the Martin kernel (see Benjamini, Pemantle and Peres [1]).

The $\geq$-half of Proposition 3.2 follows from Proposition 3.3 and the first inequality in Theorem 2.2. Specifically, on a set of diameter at most $\varepsilon$, we have $h_\varepsilon(x,y) = \log^2|x-y|/|\log\varepsilon|$, and therefore,

$$\begin{aligned}
\mathbf{P}(\mathcal{D}_\varepsilon \cap A \neq \varnothing) &\geq c\operatorname{Cap}_\varepsilon(A) &&\text{(first half of Theorem 2.2)} \\
&= c|\log\varepsilon|\operatorname{Cap}_{\log^2}(A) \\
&\asymp |\log\varepsilon|\mathbf{P}(\mathcal{I} \cap A \neq \varnothing) &&\text{(Proposition 3.3).}
\end{aligned}$$

Among the two propositions, what is left to prove is the $\leq$-half of Proposition 3.2, namely,

$$\mathbf{P}(\mathcal{D}_\varepsilon \cap A \neq \varnothing) \leq c|\log\varepsilon|\mathbf{P}(\mathcal{I} \cap A \neq \varnothing). \tag{3.13}$$

To prove this, the following corollary of Proposition 3.3 will be useful.

COROLLARY 3.4. *Let $A$ be a subset of the disk of radius $\varepsilon/2$ centered at the origin. Let $\sigma$ and $\tilde{\sigma}$ denote the respective hitting times of $B_t$ and $\tilde{B}_t$ on the circle $\{|x| = 2\varepsilon\}$. Let $z$ denote the point $(\varepsilon, 0)$ and let*

$$p = \mathbf{P}_{z,z}(A \cap B[0,\sigma] \cap \tilde{B}[0,\tilde{\sigma}] \neq \varnothing),$$
$$p' = \mathbf{P}_{z,z}(A \cap B[0,\tau_*] \cap \tilde{B}[0,\tilde{\tau}_*] \neq \varnothing)$$

*be the probabilities of two independent Brownian motions starting at $(\varepsilon, 0)$ intersecting in $A$ when stopped at $\{|x| = 2\varepsilon\}$ or $\{|x| = 3\}$ respectively. Then*

$$p' \asymp (p \cdot \log^2 \varepsilon) \wedge 1$$

*and, consequently,*

$$\mathbf{P}(\mathcal{I} \cap A \neq \varnothing) \asymp p \wedge \frac{1}{\log^2 \varepsilon}.$$

PROOF. If $|x|, |y| \leq \varepsilon/2$, then the Green function for Brownian motion stopped when it exits the disk of radius $R$ satisfying

$$G_R(x,y) \asymp \log\frac{R}{|x-y|} \tag{3.14}$$

uniformly in $R$ for $R \geq 2\varepsilon$. This follows, for instance, from $G_R(0,y) = \log(R/|y|)$ by applying a bi-Lipshitz map. Applying (3.14) to $R = 2\varepsilon$ gives

$$M_{2\varepsilon}(x,y) = \frac{G_{2\varepsilon}(x,y)}{G_{2\varepsilon}(z,y)} \asymp \frac{\log(2\varepsilon/|x-y|)}{\log(2\varepsilon/|z-y|)} \asymp \log\frac{2\varepsilon}{|x-y|}.$$

Applying (3.14) to $R=3$ then gives

$$M_3(x,y) \asymp \frac{\log(3/|x-y|)}{\log(3/|z-y|)}$$

$$\asymp \frac{\log(2\varepsilon/|x-y|) + \log(3/(2\varepsilon))}{\log(3/\varepsilon)}$$

$$\asymp 1 + \frac{M_2(x,y)}{|\log \varepsilon|}.$$

It follows that $\operatorname{Cap}_{M_3^2} \asymp 1 \wedge (\log^2 \varepsilon \cdot \operatorname{Cap}_{M_{2\varepsilon}^2})$. The first assertion of the corollary follows from this, and the second from the first and conditioning both Brownian motions to hit $D_{2\varepsilon}$. □

PROOF OF THE $\leq$-HALF OF PROPOSITION 3.2. For $z=(\varepsilon,0)$, we will show that

(3.15) $$\mathbf{P}_z(\mathcal{D}_\varepsilon \cap A \neq \varnothing) \leq cp \log^2 \varepsilon \wedge 1.$$

This suffices, since, by the Markov property,

$$\mathbf{P}(\mathcal{D}_\varepsilon \cap A \neq \varnothing) \asymp \frac{1}{|\log \varepsilon|} \mathbf{P}_z(\mathcal{D}_\varepsilon \cap A \neq \varnothing)$$

$$\leq (cp|\log \varepsilon|) \wedge \frac{1}{|\log \varepsilon|} \qquad \text{[consequence of (3.15)]}$$

$$\asymp c|\log \varepsilon| \mathbf{P}(\mathcal{I} \cap A \neq \varnothing) \qquad \text{(by Corollary 3.4)},$$

establishing (3.13) and the proposition.

To prove (3.15), let $\sigma_1 < \tau_1 < \sigma_2 < \tau_2 \cdots$ be the alternating sequence of hitting times of $\partial D_\varepsilon$ and $\partial D_{2\varepsilon}$:

$$\sigma_{n+1} = \inf\{t > \tau_n : |B_t| = \varepsilon\}$$

and

$$\tau_{n+1} = \inf\{t > \sigma_{n+1} : |B_t| = 2\varepsilon\}.$$

We call the path segments $\{B_s : \sigma_j \leq s \leq \tau_j\}$ sojourns. The left-hand side of (3.15) is bounded above by the sum

(3.16) $$\sum_{i,j \geq 1} \mathbf{P}(\sigma_i \leq \tau_*, \sigma_j \leq \tau_*, B[\sigma_i, \tau_i] \cap B[\sigma_j, \tau_j] \cap A \neq \varnothing)$$

of probabilities that sojourns $i$ and $j$ exist and intersect inside $A$. Since

$$\mathbf{P}(\tau_* < \sigma_{n+1} \mid \mathcal{F}_{\tau_n}) = \frac{\log 2}{\log(3/\varepsilon)} \asymp \frac{1}{|\log \varepsilon|},$$

the Markov property shows that the number of sojourns is geometrically distributed with mean $\log(3/\varepsilon)/\log 2 \asymp 1/|\log \varepsilon|$. For distinct sojourns, the Harnack principle again implies that the probability of their intersecting in $A$ is at most a constant multiple of $p$, and this is still true when conditioned on the number of sojourns. The expected number of pairs of sojourns is estimated by $\log^2 \varepsilon$, hence, we have a contribution of $O(p \cdot \log^2 \varepsilon)$ to the right-hand side of (3.16) from terms with $i \neq j$.

To finish, we need to estimate the probability of an $\varepsilon$-separated intersection in $A$ within a single sojourn. For $0 \leq i \leq j - 2$, let $G_{ij}$ denote the event

$$\{B_r = B_s \in S \cap A \text{ for some } r \in [i\varepsilon^2/2, (i+1)\varepsilon^2/2] \text{ and } s \in [j\varepsilon^2/2, (j+1)\varepsilon^2/2]\}.$$

Let $t_j = (j - \frac{1}{2})\frac{\varepsilon^2}{2}$ and let $\sigma$ be the hitting time of $\{|x| = 2\varepsilon\}$. We apply the Markov property at time $\sigma_n$ to estimate the summand in (3.16) with $i = j = n$, then sum over $n$. This bounds the contributions to (3.16) from off-diagonal terms by

$$\sum_n \mathbf{P}_z(\sigma_n < \tau_*) \sum_{i<j} \mathbf{P}_{B(\sigma_n)}(G_{ij}, t_j < \sigma).$$

The sum over $n$ is $O(|\log \varepsilon|)$ and the sum over $0 \leq i < j$ of $\mathbf{P}_{B(\sigma_n)}(t_j < \sigma)$ is $O(1)$ (e.g., this is at most the sum of $j$ times the probability that $\{|x| = 2\varepsilon\}$ is not hit by time $j\varepsilon^2/2$, which is at most the expected square of the time for a Brownian motion to reach $\{|x| = 4\}$). We will be done, therefore, when we have shown that

$$(3.17) \qquad \sup_{0 \leq i \leq j-2, |z|=\varepsilon} \mathbf{P}_z(G_{ij} \mid \mathcal{F}_{t_j}) \leq cp$$

on the event $\{t_j < \sigma\}$ (actually, an upper bound of $cp|\log \varepsilon|$ would suffice).

This is more or less obvious from the Markov property, but we go ahead and spell out the details. Let $\omega_i$ denote the $i$th sub-sojourn defined by $\omega_i(s) = \omega(s + i\varepsilon^2/2)$ for $0 \leq s \leq \varepsilon^2/2$. Let $\mu_{ij}$ denote the conditional law of $\omega_i$ under $\mathbf{P}_z$ given $\mathcal{F}_{t_j}$ and $\mu$ denote the $P_z$-law of $\omega$ on the interval $[\varepsilon^2/2, \varepsilon^2]$. The quantity $p$ is estimated by the probability of two independent draws from $\mu$ intersecting inside $A$; conditioning on $\mathcal{F}_{t_j}$ makes $\omega_i$ and $\omega_j$ independent, so (3.17) follows if we can show that

$$(3.18) \qquad \frac{d\mu_{ij}}{d\mu} \leq C \mathbf{1}_{t_j < \sigma}$$

when $0 \leq i \leq j - 2$ or $i = j$. For $i = j$, $\mu_{jj}$ and $\mu$ are Wiener measure from starting points with comparable densities. For $1 \leq i \leq j - 2$, we use the Markov property to write

$$\mu_{ij} = \int \mu^{xy}_{\varepsilon^2/2} \mathbf{1}_H \, d\pi_{ij}(x, y),$$

$$\mu = \int \mu^{xy}_{\varepsilon^2/2} \, d\pi(x, y),$$

where $\mu_t^{xy}$ is the law of a Brownian bridge from $x$ to $y$ in time $t$, $H$ is the event that the path remains inside the ball of radius $2\varepsilon$, and $\pi_{ij}$ and $\pi$ are mixing measures. By Bayes' rule and the Markov property,

$$\frac{\pi_{i,j}(x,y)}{\pi(x,y)} = \frac{1}{Z}\,\mu_{x,y}(H)\mu_{i\varepsilon^2/2}^{zx}(H)\mu_{(t_j-i-1)\varepsilon^2/2}^{y,B(t_j)}(H),$$

where $Z$ is the normalizing constant gotten by integrating the product of the three probabilities on the right-hand side against $\pi(x,y)$. The probabilities are all at most 1, so all we need is that $Z$ is at least $c > 0$. By Brownian scaling, we see that the three probabilities are at least a constant when $|x|,|y| < \varepsilon$ and, since $\pi$ gives positive measure to this set, the verification for $1 \le i \le j-2$ is complete. Finally, for $i = 0$, we compare to $\mu'$ instead of $\mu$, where $\mu'$ is the $\mathbf{P}_z$-law of $\omega$ on $[0,\varepsilon^2/2]$. This establishes (3.17) and, hence, (3.15) and the remainder of Proposition 3.2. □

PROOF OF THE SECOND INEQUALITY IN THEOREM 2.2.  Let

$$\tau = \tau_\varepsilon = \inf\{t : B_t \in A \text{ and } B_s = B_t \text{ for some } s \le t - \varepsilon^2\}$$

be the first time that a point of $A$ is hit by the Brownian motion and has previously been hit at a time at least $\varepsilon^2$ in the past; thus, $\mathbf{P}(\tau \le \tau_*) = \mathbf{P}(\mathcal{D}_\varepsilon \cap A \ne \varnothing)$. The second inequality in Theorem 2.2 is equivalent to the existence of a measure $\nu$ on $A$ whose mass is equal to $\mathbf{P}(\tau < \tau_*)$ and whose energy is at most a constant multiple of this [normalizing $\nu$ to be a probability measure gives an energy of $C\mathbf{P}(\mathcal{D}_\varepsilon \cap A \ne \varnothing)^{-1}$, thereby witnessing the inequality].

To construct $\nu$, partition the plane into a grid of squares of side $\varepsilon/3$. For each square $S$ in the grid, let $\nu_S$ be a probability measure of minimal $\log^2$-energy on $S \cap A$. By Proposition 3.3,

$$\mathcal{E}_{h_\varepsilon}(\nu_S) = \frac{1}{|\log \varepsilon|}\mathcal{E}_{\log^2}(\nu_S) \le C\mathbf{P}(\mathcal{I} \cap S \cap A \ne \varnothing)^{-1}$$

and so by Proposition 3.2, for a different constant,

(3.19) $$\mathcal{E}_{h_\varepsilon}(\nu_S) \le C\mathbf{P}(\mathcal{D}_\varepsilon \cap S \cap A \ne \varnothing)^{-1}.$$

Let

$$\nu := \sum_S \mathbf{P}(B_\tau \in S, \tau < \tau_*)\nu_S.$$

Clearly, we have constructed $\nu$ so that $\|\nu\| = \mathbf{P}(\tau < \tau_*)$. It remains to show that $\mathcal{E}_{h_\varepsilon}(\nu) \le c\mathbf{P}(\tau < \tau_*)$. We will tally separately the contributions to the energy from pairs $(x,y)$ at distances at least $\varepsilon$ and at most $\varepsilon$, showing

(3.20) $$\int h_\varepsilon(x,y)\mathbf{1}_{|x-y|\ge\varepsilon}\,d\nu(x)\,d\nu(y) \le C\|\nu\|$$

and

$$\int h_\varepsilon(x,y)\mathbf{1}_{|x-y|\leq\varepsilon}\,d\nu(x)\,d\nu(y) \leq C\|\nu\|. \tag{3.21}$$

For the bound (3.20) on the first piece, observe that points separated by $\varepsilon$ are in nonadjacent squares $S$ and $S'$, and that the value of $h$ at any $x \in S$ and $y \in S'$ is estimated by the $|\log|x_* - y_*||$ for any $x_* \in S, y_* \in S'$. Therefore, on the event $\{B_\tau \in S\}$, we may replace $x \in S$ by $B_\tau$ to obtain

$$\int h_\varepsilon(x,y)\mathbf{1}_{|x-y|\geq\varepsilon}\,d\nu(x)\,d\nu(y)$$
$$\leq C \sum_{S'} \mathbf{P}(B_\tau \in S', \tau < \tau_*)$$
$$\times \int d\nu_{S'}(y) \left[ \sum_{S \text{ not adjacent to } S'} \mathbf{E}(\mathbf{1}_{B_\tau \in S, \tau < \tau_*} |\log|B_\tau - y||) \right]$$
$$\leq \|\nu\| \sup_y V(y),$$

where

$$V(y) = \mathbf{E}(|\log|B_\tau - y||\mathbf{1}_G)$$

is the logarithmic potential at $y$ of the subprobability law of $B_\tau$ restricted to the event $G := \{\tau < \tau_*, |B_\tau - y| \geq \varepsilon)\}$.

To see that $V(y)$ is bounded, fix $y$ and observe that the probability that $\mathcal{D}_\varepsilon$ intersects the $\delta$-ball $D_y$ is at least equal to the probability that it does so after time $\tau$ has been reached. Throwing away those paths where $B_\tau$ is within $\varepsilon$ of $y$, we have, by the Markov property and (3.2),

$$\mathbf{P}(\mathcal{D}_\varepsilon \cap D_y \neq \varnothing) \geq c\mathbf{E}\frac{\log\varepsilon \log|X - y|}{\log^2 \delta}\mathbf{1}_G.$$

On the other hand, by (3.1),

$$\mathbf{P}(\mathcal{D}_\varepsilon \cap D_y \neq \varnothing) \asymp \frac{|\log\varepsilon|}{\log^2 \delta}.$$

It follows that $\mathbf{E}|\log|X - y||\mathbf{1}_G \leq c^{-1}$, which is the desired bound on the first piece.

For the bound (3.21) on the second piece, begin with the well known trick of reducing to the diagonal:

$$\int h_\varepsilon(x,y)\mathbf{1}_{|x-y|\leq\varepsilon}\,d\nu(x)\,d\nu(y) \leq C\sum_S \mathbf{P}(B_\tau \in S, \tau < \tau_*)^2 \mathcal{E}_{h_\varepsilon}(\nu_S). \tag{3.22}$$

One way to see this is to observe that, while $|x - y| \leq \varepsilon$ and $x \in S$ does not force $y \in S$, it does force $y$ to be in one of 49 nearby squares. The function $\log^2 |x-y|/|\log \varepsilon|$ is positive definite, so one may use the Cauchy–Schwarz inequality to conclude (3.22). In fact, (3.22) holds when $h_\varepsilon$ is not positive definite but only assumed to be monotone; see Pemantle and Peres [4], equation 11 for details.

Finally, since $\mathbf{P}(B_\tau \in S, \tau < \tau_*) \leq \mathbf{P}(\mathcal{D}_\varepsilon \cap S \cap A \neq \varnothing)$, we see from (3.19) that

$$\mathbf{P}(B_\tau \in S, \tau < \tau_*)^2 \mathcal{E}_{h_\varepsilon}(\nu_S) \leq C \mathbf{P}(B_\tau \in S, \tau < \tau_*).$$

Summing over $S$ bounds the right-hand side of (3.22) by $\mathbf{P}(\tau < \tau_*) = \|\nu\|$, establishing (3.21) and finishing the proof of Theorem 2.2. □

**4. Proof of Theorem 2.4.** There are two obvious choices for the set $A_2$. The first is the set $\mathcal{P}$ of points $x$ such that a Brownian motion started at $x$ and run for any positive time almost surely has a double point in $A$. Call such a point an *immediate point*. The second choice would be the set $\mathcal{R}$ of regular points of $A$ with respect to the potential of the least-energy measure for the kernel $K(x,y) = \log^2 |x-y|$. [A *regular point* $x$ for the potential $\int K(x,y) \, d\nu(y)$ of a measure $\nu$ is one where the potential reaches its maximum value.] If $\mathcal{P} = \mathcal{R}$, then Theorem 2.4 has a very short proof:

> Let $A_2 = \mathcal{P} = \mathcal{R}$. It is well known (see Proposition 4.3 below) that the non-regular points $A_1 := A \setminus \mathcal{R}$ must have zero $K$-capacity, and thus, using the intersection criterion from Fitzsimmons and Salisbury [3], cannot intersect $\mathcal{D}$. This is property (1) required by the theorem. But property (2) in the Theorem is satisfied by definition of $\mathcal{P}$, noting that by what we just proved, having a double point in $A$ is the same as having a double point in $A_2$.

Embarrassingly, we do not know whether $\mathcal{R} = \mathcal{P}$. We can, however, establish something close, namely, Lemma 4.1, which will be enough to prove the theorem. The apparent obstacle to proving the equality of $\mathcal{P}$ and $\mathcal{R}$ is their different nature: $\mathcal{P}$ is defined probabilistically and the definition is inherently local, while $\mathcal{R}$ is defined analytically and its definition is at first glance nonlocal. Accordingly, we define an analytic version of $\mathcal{P}$ and a localized version of $\mathcal{R}$ as follows.

Fix the closed set $A$ and let $\xi$ be a point of $A$. Let $f$ be any decreasing continuous function from $\mathbb{R}^+$ to $\mathbb{R}^+$ going to infinity at 0, and let $M_\xi$ denote the $f$-Martin kernel at $\xi$:

$$M_\xi(x,y) := \frac{f(|x-y|)}{f(|\xi-y|)}.$$

We say that $A$ has nonvanishing local Martin capacity (NLMC) at $\xi$ if and only if

$$\lim_{\varepsilon \to 0} \operatorname{Cap}_{M_\xi}(A \cap \{y : |y - \xi| < \varepsilon\}) > 0.$$

Let $\mathcal{P}'$ denote the set of points with NLMC. The relation to $\mathcal{P}$ will be clarified shortly.

Call a point $\xi \in A$ *strongly regular* if and only if the $f$-capacity of $A \cap \{y : |y - \xi| < \varepsilon\}$ is nonzero for every $\varepsilon$, and $\xi$ is a regular point for the potential of the least $f$-energy measure on each such set. Let $\mathcal{R}'$ denote the set of strongly regular points.

LEMMA 4.1 (Strongly regular implies NLMC for any gauge). *For any $A$ and $f$ as above, the inclusion $\mathcal{R}' \subseteq \mathcal{P}'$ holds.*

PROOF. Fix $\xi \in \mathcal{R}'$. Given any ball $D$ containing $\xi$, let $\nu_D$ denote the measure minimizing the $f$-energy and let $\Phi_D$ denote its potential:
$$\Phi_D(x) = \int f(|x - y|) \, d\nu_D(y).$$
By assumption, $\Phi_D(\xi)$ is equal to the maximum value of $\Phi_D$. It is well known that the maximum value is attained on a set of full measure; standard references such as Carleson [2] state unnecessary assumptions on $f$, so we include the proof (Proposition 4.3 below). It follows that
$$\mathcal{E}_f(\nu_D) = \Phi_D(\xi).$$
Define a new measure $\rho_D$, which is a probability measure, by
$$\frac{d\rho_D}{d\nu_D}(y) = \frac{f(|\xi - y|)}{\Phi_D(\xi)}.$$
The potential of this new measure with respect to the Martin kernel $M_\xi$ at a point $x$ is computed to be
$$\frac{1}{\Phi_D(\xi)} \int M_\xi(x, y) f(|\xi - y|) \, d\nu_D(y) = \frac{1}{\Phi_D(\xi)} \int f(|x - y|) \, d\nu_D(y) = \frac{\Phi_D(x)}{\Phi_D(\xi)}.$$
Since $\xi$ is regular for $\Phi_D$, this is at most 1. Since the Martin potential is bounded by 1, the Martin energy $\mathcal{E}_{M_\xi}(\nu_D)$ of the probability measure $\nu_D$ is also at most 1, and we see that each ball $D$ has Martin capacity at least 1. □

LEMMA 4.2 (NLMC points are immediate). *Let the closed set $A$ have nonvanishing local Martin capacity at $\xi$ for the $\log^2$ Martin gauge*
$$M_\xi(x, y) := \log^2 |x - y| / \log^2 |\xi - y|.$$
*Then $\xi$ is an immediate point.*

REMARK. We first remark that if $\Xi$ is the range of a transient Markov process with Green function $G$ and $M_\xi$ is the Martin kernel for the process started at $\xi$, then the implication holds in both directions: the set $A$ has nonvanishing local $M_\xi$-capacity near $\xi$ if and only if the process started from $\xi$ almost surely intersects $A$ in any positive time interval. This follows from the methods of Benjamini, Pemantle and Peres [1].

The set of double points is not the range of a Markov process, which makes proving a reverse implication tricky, but the direction in the lemma may still be obtained by applying the method of second moments. Recall that $H(\Lambda, \varepsilon)$ denotes the event that there is a double point in the set $\Lambda$ with an $\varepsilon^2$ time separation.

PROOF OF LEMMA 4.2. Begin by observing it is enough to show

$$(4.1) \qquad \mathbf{P}_\xi(\mathcal{D} \cap A \neq \varnothing) \geq c \operatorname{Cap}_{M_\xi}(A).$$

For, under the hypothesis of NLMC, this implies that

$$\inf_{\varepsilon > 0} \mathbf{P}_\xi[H(A \cap \{y : |y - \xi| < \varepsilon\}, \varepsilon)] > 0.$$

By Fatou's lemma,

$$\mathbf{P}_\xi \left[ \limsup_{\varepsilon \to 0} H(A \cap \{y : |y - \xi| < \varepsilon\}, \varepsilon) \right] > 0,$$

whence, with positive probability, $\mathcal{D}$ intersects $A$ in a set with $\xi$ as a limit point. Since

$$\mathbf{P}_\xi(|B_t - \xi| < \varepsilon \text{ for some } t > s) \to 0$$

as $\varepsilon \to 0$ for any fixed $s$, it follows that a Brownian motion run from $\xi$ for an arbitrarily short time has a double point in $A$ with probability bounded away from zero. By Blumenthal's zero–one law, this probability must be 1, so $\xi$ is an immediate point.

We will prove something slightly stronger than (4.1), replacing $\mathcal{D}$ in (4.1) by a subset akin to $\mathcal{D}_\varepsilon$ but where the value of $\varepsilon$ depends on the distance to the point $\xi$:

$$\mathcal{D}_* = \{x : B_s = B_t \text{ for some } s, t < \tau \text{ with } |x - \xi|^2 \leq t - s \leq |x - \xi|\}.$$

(Recall we stop at the time $\tau$ that the Brownian motion exits a disk of radius 3.) Let $H_*(S)$ denote the event that $\mathcal{D}_*$ has nonempty intersection with $S$, and let $S_x$ denote the disk of radius $\delta |x - \xi|$ centered at $x$. The relevant two-point correlation estimate we will prove is, for $|x - \xi| \leq |y - \xi|$,

$$(4.2) \qquad \frac{\mathbf{P}_\xi[H_*(S_x) \cap H_*(S_y)]}{\mathbf{P}_\xi(H_*(S_x)) \mathbf{P}_\xi(H_*(S_y))} \leq C M_\xi(x, y).$$

Assuming this, the proof is finished in the same manner as the proof of the lower bound in Theorem 2.2, as follows.

Let $\mu$ be any probability measure on $A$. Fix $1/4 > \delta > 0$, which will later be sent to zero. According to (1.2), we may assume $A$ to be a finite disjoint union of squares of a lattice which has been subdivided so that squares at distance $r$ from $\xi$ have sides between $\delta r$ and $3\delta r$; the Whitney decomposition of the complement of $\xi$ forms such a subdivision. This contains the union of disks $\{S_x : x \in \mathcal{B}\}$ and, as before, we may choose $\mathcal{B}$ so no two disks are closer to each other than the radius of the smaller disk, while the union of the disks still has measure at least $c\mu(A)$. Define

$$X := \sum_{x \in \mathcal{B}} \frac{1}{\mathbf{P}(H_*(S_x))} \mu(S_x) \mathbf{1}_{H_*(S_x)}.$$

Then $\mathbf{E}X \geq c$ and by (4.2),

$$\mathbf{E}X^2 \leq 2C \sum_{S_x, S_y} \mu(S_x) \mu(S_y) M_\xi(x, y).$$

Here, instead of counting each pair twice, we have summed over $(x, y)$ for which $|x - \xi| \leq |y - \xi|$ and then doubled. As in (3.12), for $x' \in S_x$ and $y' \in S_y$, we have $M_\xi(x', y') \asymp M_\xi(x, y)$, so we may apply the second moment method to obtain

$$\mathbf{P}(H_*(A)) \geq \frac{(\mathbf{E}X)^2}{\mathbf{E}X^2} \geq c^2(c + 2\mathcal{E}(\mu))^{-1}.$$

This is uniform in $\delta$, so sending $\delta$ to zero proves (4.1). It remains to prove (4.2).

Given $x$ and $y$ and $\delta \leq 1/4$, observe that when $|x - \xi| < |y - \xi|^2$, then $\mathbf{P}_\xi$ makes $H_*(S_x)$ and $H_*(S_y)$ independent up to a constant factor which is independent of $\delta$. To see this, compute the probabilities of hitting in various orders to find that the dominant term comes from hitting $S_x$ twice before the Brownian motion reaches a disk of radius $|y - \xi|/2$; after this, the conditional probability of $H_*(S_y)$ is only a constant multiple of the unconditional probability. Independence up to a constant factor means a two-point correlation function bounded by a constant, whence (4.2) is satisfied.

In the complementary case, the ratio of $\log|x - \xi|$ to $\log|y - \xi|$ is bounded, so we may again compute the two-point correlation function as in the proof of Theorem 2.2. Recall from (3.2) of Lemma 3.1 that using $\mathbf{P}_\xi$ instead of $\mathbf{P}$ boosts the individual probabilities of $H(D_x)$ by a factor of $|\log|x - \xi||$. The same holds for $H_*(S_x)$. Thus,

$$\mathbf{P}_\xi(H_*(S_x)) \asymp \frac{\log^2|x - \xi|}{\log^2(\delta|x - \xi|)}.$$

The probability of $H_*(S_x) \cap H_*(S_y)$ is again computed by summing the probabilities of various scenarios, the likeliest of which (up to a constant factor) is a hit on $S_x$, then on $S_y$, then a time separation of at least $|x - \xi|^2$, then another hit on $S_x$ and then on $S_y$. Multiplying this out gives

$$\frac{\log|x - \xi|}{\log(\delta|x - \xi|)} \cdot \frac{\log|x - y|}{\log(\delta|y - \xi|)} \cdot \frac{\log|x - \xi|}{\log(\delta|x - \xi|)} \cdot \frac{\log|x - y|}{\log(\delta|y - \xi|)},$$

PROBABILITY OF INTERSECTING BROWNIAN DOUBLE POINTS     19

which results in the estimate (4.2).  □

For completeness' sake, as mentioned above, we repeat here the standard argument to show that the complement of the strongly regular points is a set of zero capacity.

PROPOSITION 4.3 ($\mathcal{R}^c$ has zero capacity in any gauge). *The set $\mathcal{R}^c$ of nonregular points of a set $A$ for the minimizing measure with respect to any continuous gauge $f$ has zero $f$-capacity (and, in particular, has zero minimizing measure). It follows from countable additivity that $\operatorname{Cap}_f(\mathcal{R}')^c = 0$ as well.*

PROOF. Assume to the contrary that $A \setminus \mathcal{R}$ has positive capacity. Let $\nu$ be a minimizing probability measure on $A$ for $\mathcal{E}_f$. Then for some $\delta$, the set $\{y \in A : \Phi_\nu(y) < (1-\delta)\mathcal{E}(\nu)\}$ has positive capacity, where $\Phi_\nu(y) := \int f(x,y)\,d\nu(x)$ is the $f$-potential of $\nu$ at $y$. Fix such a $\delta$ and let $\mu$ be a probability measure supported on this set with $\mathcal{E}_f(\mu) < \infty$. For $\varepsilon \in (0,1)$, consider the measure $\rho_\varepsilon := (1-\varepsilon)\nu + \varepsilon\mu$. Its energy is given by

$$(1-\varepsilon)^2 \mathcal{E}_f(\nu) + \varepsilon^2 \mathcal{E}_f(\mu) + 2\varepsilon(1-\varepsilon) \int\int f(x,y)\,d\mu(x)\,d\nu(y).$$

The double integral is equal to $\int \Phi_\nu(x)\,d\mu(x)$ and since this is at most $(1-\delta)\mathcal{E}_f(\nu)$ on the support of $\mu$, the energy of $\rho_\varepsilon$ is bounded above by

$$[(1-\varepsilon)^2 + 2\varepsilon(1-\varepsilon)(1-\delta)]\mathcal{E}_f(\nu) + \varepsilon^2 \mathcal{E}_f(\mu).$$

Write this as $\mathcal{E}_f(\nu)(1 - 2\varepsilon\delta + \varepsilon^2 Q)$, where $Q = \mathcal{E}_f(\mu)/\mathcal{E}_f(\nu) + 2\delta - 1 < \infty$, and take the derivative at $\varepsilon = 0$ to see that $\mathcal{E}_f(\rho_\varepsilon) < \mathcal{E}_f(\nu)$ for small positive $\varepsilon$. This contradicts the minimality of $\mathcal{E}_f(\nu)$ and proves the proposition.  □

Finally, we complete the proof of the decomposition as follows. Let $A_2$ be the set of strongly regular points of $A$. We have just seen that $A_1 := A \setminus A_2$ has zero capacity in the gauge $\log^2 |x-y|$. By Fitzsimmons and Salisbury [3], this implies that $A_1$ is almost surely disjoint from the set of Brownian double points, which is property (1).

On the other hand, by Lemma 4.1, $A$ has NLMC at each point of $A_2$, and by Lemma 4.2, all such points are immediate for $A$. Using the fact that $A_1$ has no double points again, we conclude that property (2) in the statement of Theorem 2.4 is satisfied.

DEPARTMENT OF MATHEMATICS
UNIVERSITY OF PENNSYLVANIA
209 S. 33RD STREET
PHILADELPHIA, PENNSYLVANIA 19104
USA
E-MAIL: pemantle@math.upenn.edu

DEPARTMENT OF STATISTICS
UNIVERSITY OF CALIFORNIA
EVANS HALL
BERKELEY, CALIFORNIA 94720
USA
E-MAIL: peres@stat.berkeley.edu